\documentclass{article}

\usepackage{authblk}
\usepackage{amsmath, amssymb, amsthm,bbm}
\usepackage{graphicx,url,enumerate,xcolor,eucal}

\newtheorem{theorem}{Theorem}[section]
\newtheorem{lemma}[theorem]{Lemma}
\newtheorem{proposition}[theorem]{Proposition}
\newtheorem{corollary}[theorem]{Corollary}
\theoremstyle{definition}
\newtheorem{definition}[theorem]{Definition}
\newtheorem{example}[theorem]{Example}

\newcommand{\zexp}{\mathfrak{exp}\mspace{1mu}}

\newcommand{\zco}{\mathfrak{C}\mspace{1mu}}

\newcommand{\CZ}{\mathbb{C}\mathfrak{Z}}

\newcommand{\Wexp}{\Omega\mathfrak{exp}\mspace{1mu}}

\DeclareMathOperator{\tr}{tr}

\DeclareMathOperator{\per}{per}

\newcommand{\zm}{\mathbf{0}}
\newcommand{\mI}{\mathbb{I}}

\numberwithin{equation}{section}

\title{Bipartite Matchings, Matrix Permanents, and Zeon Power Series}

\author{Muhammad Waqar Ul Hassan}
\author[1*]{G. Stacey Staples
\affil[*]{Corresponding author: \texttt{sstaple@siue.edu}}}
\affil[1]{Department of Mathematics and Statistics,
Southern Illinois University Edwardsville,
Edwardsville, Illinois, USA}

\begin{document}

\maketitle

\begin{abstract}
Perfect matchings of a bipartite graph are counted by the permanent of the graph's biadjacency matrix, and permanents of sequences of structured biadjacency matrices lead to integer sequences.  We consider a number of such sequences and illustrate how zeon-algebraic representations of bipartite graphs lead to recursively defined zeon sequences.  Then, using formal power series methods, we show that  the ``$\Omega$-exponentials'' of these zeon sequences reveal numbers of perfect matchings in the graphs.  
\end{abstract}

\section{Introduction}\label{intro}

The permanent is a scalar-valued function defined for square matrices and is
closely related to the determinant. The principal difference between them lies
in the treatment of permutations: each term in the determinant is multiplied
by the sign of its corresponding permutation, whereas the permanent omits
this sign factor \cite{marcus1965}. Computing a permanent is generally more
difficult than computing a determinant \cite{brualdi1991}. However, for
certain structured matrices, including particular classes of tridiagonal
matrices, explicit formulas may simplify the computation. Finding exact
formulas and efficient methods for computing the permanents of structured
matrices therefore remains an important problem. In this paper, the results
are established using induction, Laplace expansion, recurrence relations, and
known permanent identities.

The study of permanents also extends to several specialized matrix classes
and applications. Armario adapted Ryser's formula to the structure of
Sylvester Hadamard matrices, reducing the number of subsets of distinct rows
that must be considered in the computation \cite{armario2021}. Singh and
Bapat developed recurrence relations based on Laplace expansion for computing
permanents and permanent polynomials using induced subdigraphs associated with
square matrices \cite{singh2017}. Permanent-based graph models have also been
applied to decision-making problems \cite{baykasoglu2014}.

In this paper, we consider families of bipartite graphs whose perfect matchings are counted by integer sequences and illustrate how zeon algebra representations of bipartite graphs lead to recursively defined zeon sequences.  Using formal power series methods, we show that  the ``$\Omega$-exponentials'' of these zeon sequences reveal numbers of perfect matchings in the graphs. 

The rest of Section \ref{intro} contains essential terminology, notation, and fundamental results on bipartite graphs, biadjacency matrices, integer sequences, and counting matchings via the permanent of the biadjacency matrix.  Section \ref{matchings and matrices} contains results concerning the permanents of some specific families of biadjacency matrices and the numbers of perfect matchings in their corresponding bipartite graphs.   

Section \ref{zeon formulations} introduces zeon algebras\footnote{Zeon algebras can be interpreted as commutative fermion algebras.  Equivalently, they can be constructed as subalgebras of Clifford algebras.} and provides matrix-free representations of bipartite graphs that allow perfect matchings to be counted by considering powers and exponentials of zeon sequences.  In Section \ref{OmegaExp}, we introduce the ``$\Omega$-exponential'' of a zeon sequence $\Gamma=(\gamma_j)$ and show how coefficients of $\Wexp(t\Gamma$) can be used to count the perfect matchings of a bipartite graph.  As examples, we use zeon methods to recover the results established in Section \ref{matchings and matrices}.

The paper concludes with Section \ref{conclusion}, where we summarize our results and suggest some avenues for further investigation.

\begin{description}
\item [Notation]
In this paper, we consider four types of sequences: integers, matrices, graphs, and zeon elements.  In order to avoid confusion, we adopt the following notational conventions. \begin{itemize}
\item Integer sequences will be denoted using lower case Fraktur characters; e.g., $(\mathfrak{f}_n)$ will denote the Fibonacci sequence.
\item Matrix sequences will be denoted with capital Roman letters; e.g., $(M_n)$.
\item Graph sequences will be indicated with capital Roman letters in caligraphic font; e.g., $(\mathcal{G}_n)$.
\item Zeon algebra sequences will be indicated with lower case Greek characters; e.g., $(\alpha_n)$.
\end{itemize}
\end{description}

\subsection{Permanent of a Square Matrix}\label{permanent intro}

\begin{definition}
The {\em permanent} of an $n\times n$ matrix $Z=(z_{ij})$ is defined by
\begin{equation*}
\per(Z)
=\sum_{\sigma\in S_n}\prod_{i=1}^{n}z_{i,\sigma(i)},
\end{equation*}
where $S_n$ is the symmetric group of degree $n$.
\end{definition}

For example, if
$Z=\begin{pmatrix}p&q\\ r&s\end{pmatrix}$, then
$\per(Z)=ps+qr$.

The permanent of an $n\times n$ matrix $A=(a_{ij})$ can be computed
recursively using a Laplace expansion \cite{marcus1965}. Unlike the
corresponding determinant expansion, the permanent expansion does not
contain the sign factor $(-1)^{i+j}$. Expanding along the first row gives
\begin{equation*}
\per(A)
=\sum_{j=1}^{n}a_{1j}\per(A(1\vert j)),
\end{equation*}
where $A(1\vert j)$ is the $(n-1)\times(n-1)$ submatrix obtained by
deleting the first row and the $j$-th column of $A$. This process is
repeated recursively until a $1\times1$ matrix is reached, for which
\begin{equation*}
\per((a_{11}))=a_{11}.
\end{equation*}

\subsubsection{Permanents of Matrix Sequences}

We are particularly interested in sequences of matrices representing bipartite graphs of increasing size.  In particular families of sequences, the permanents satisfy recurrence relations that simplify the permanent computations.  As an example, we consider tridiagonal matrices.

Let $\{Z_n\}_{n\geq1}$ be a sequence of $n\times n$ tridiagonal matrices of the form
\begin{equation*}
Z_n=\begin{pmatrix}
z_{11} & z_{12} & 0 & \cdots & 0\\
z_{21} & z_{22} & z_{23} & \ddots & \vdots\\
0 & z_{32} & z_{33} & \ddots & 0\\
\vdots & \ddots & \ddots & \ddots & z_{n-1,n}\\
0 & \cdots & 0 & z_{n,n-1} & z_{nn}
\end{pmatrix}.
\end{equation*}

As seen in \cite{dafonseca2010}, the permanents of these matrices satisfy the  following recurrence: $\per(Z_1)=z_{11}$, $\per(Z_2)=z_{11}z_{22}+z_{12}z_{21}$, and
\begin{equation*}
\per(Z_n)=z_{nn}\per(Z_{n-1})+z_{n-1,n}z_{n,n-1}\per(Z_{n-2})
\end{equation*}
for $n\ge 3$.

\begin{example}
Consider the sequence $(H_n)$ of $n\times n$ matrices whose $n$th term $H_n=(h_{ij})$ is defined by \begin{equation*}
h_{ij}=
\begin{cases}
1, & \text{if }|i-j|\leq1,\\
0, & \text{otherwise}.
\end{cases}
\end{equation*}
Thus,
\begin{equation*}
H_n=
\begin{pmatrix}
1 & 1 & 0 & \cdots & 0\\
1 & 1 & 1 & \ddots & \vdots\\
0 & 1 & 1 & \ddots & 0\\
\vdots & \ddots & \ddots & \ddots & 1\\
0 & \cdots & 0 & 1 & 1
\end{pmatrix}.
\end{equation*}
The permanents satisfy $\per(H_1)=1$, $\per(H_2)=2$, and 
\begin{equation*}
\per(H_n)
=\per(H_{n-1})+\per(H_{n-2}),
\qquad n\geq3.
\end{equation*}
Hence, the permanents of matrices of larger order can be computed
efficiently using this recurrence.
\end{example}

Going forward, we will consider sequences of square matrices representing bipartite graphs and derive recurrence relations for their permanents, revealing numbers of perfect matchings in the graphs. 

\subsection{Biadjacency Matrices and Perfect Matchings of Bipartite Graphs}\label{bipartite graphs subsection}

\begin{figure}[h]
    $$\includegraphics[height=100pt]
    {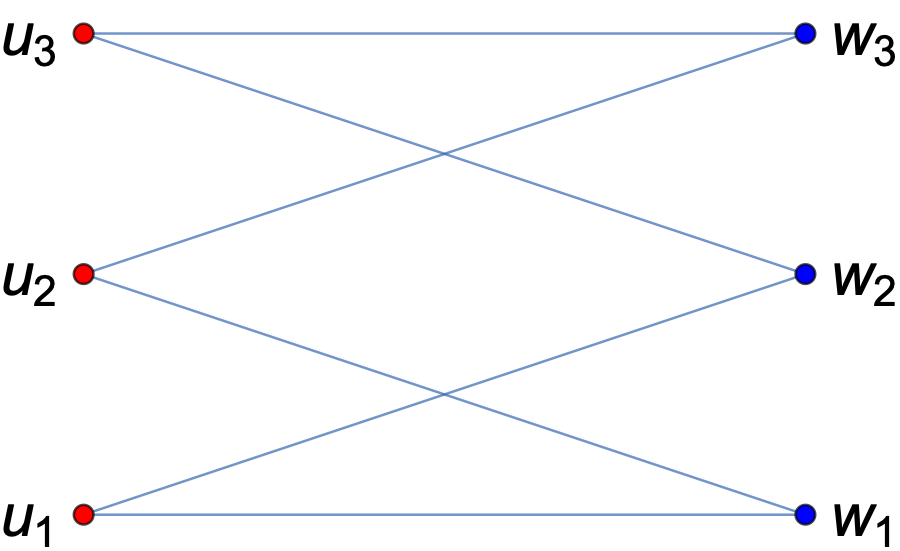}$$
    \caption{A bipartite graph with bipartition $(U,W)$.}
    \label{fig:bipartite}
\end{figure}

Formally, a graph $G=(V,E)$ is {\em bipartite} if its vertex set $V$ can be partitioned into two disjoint nonempty subsets $U$ and $W$ such that every edge in $E$ has one endpoint in $U$ and the other in $W$. Thus,
\begin{equation*}
V=U\cup W
\qquad \text{and} \qquad
U\cap W=\varnothing.
\end{equation*}
The pair $(U,W)$ is called a \emph{bipartition} of $G$.  We denote the disjoint union of $U$ and $W$ by $U\sqcup W$.  A simple example can be seen in Figure \ref{fig:bipartite}, where $U=\{u_1, u_2, u_3\}$ and $W=\{w_1, w_2, w_3\}$.

Let $G=(V,E)$ be a bipartite graph with bipartition $(U,W)$, where
$U=\{u_1,u_2,\ldots,u_m\}$ and $W=\{w_1,w_2,\ldots,w_n\}$. The
{\em biadjacency matrix} of $G$ is the $m\times n$ matrix
$B(G)=(b_{ij})$ defined by
\begin{equation*}
b_{ij}=
\begin{cases}
1, & \text{if }\{u_i,w_j\}\in E,\\
0, & \text{otherwise}.
\end{cases}
\end{equation*}

If the vertices are ordered as
$(u_1,\ldots,u_m,w_1,\ldots,w_n)$, then the full adjacency matrix of
$G$ is
\begin{align*}
A(G)&=
\begin{pmatrix}
\zm_{m\times m} & B(G)\\
B(G)^\intercal & \zm_{n\times n}
\end{pmatrix}.
\end{align*}

A {\em matching} in a graph $G=(V,E)$ is a subset of edges of $G$ in which no pair is coincident with a common vertex.  A {\em perfect matching} of a graph is a matching of cardinality $\lfloor |V|/2\rfloor$.  In particular, when $G$ is a graph on $2n$ vertices, a perfect matching has cardinality $n$.  

A {\em balanced bipartite graph} is a bipartite graph $G=(V,E)$ whose bipartition $V=U\sqcup W$ satisfies $|U|=|W|$.   A balanced bipartite graph's perfect matchings are counted by the biadjacency matrix 
permanent as follows \cite{oteles2017}.

\begin{proposition}\label{prop:perfect-matching}
Let $G=(V,E)$ be a bipartite graph with bipartition $(U,W)$, where
$|U|=|W|=n$. If $\mu(G)$ denotes the number of perfect matchings of
$G$, then
\begin{equation*}
\mu(G)=\per(B(G)).
\end{equation*}
\end{proposition}

\begin{proof}
For each $\sigma\in S_n$, the term
$\prod_{i=1}^{n}b_{i,\sigma(i)}$ equals $1$ precisely when the edges
$\{u_i w_{\sigma(i)}:1\leq i\leq n\}$ form a perfect matching of $G$. Thus, the nonzero terms in $\per(B(G))$ are in one-to-one correspondence with the perfect matchings of $G$.
\end{proof}

In this case, 
\begin{align*}
A(G)&=
\begin{pmatrix}
\zm_{n\times n} & B(G)\\
B(G)^\intercal & \zm_{n\times n}
\end{pmatrix}\\
&=\begin{pmatrix}0&1\\0&0\end{pmatrix}\otimes B(G)+\begin{pmatrix}0&0\\1&0\end{pmatrix}\otimes B(G)^\intercal.
\end{align*}

Given the structure of the adjacency matrix, the following corollary is immediate.

\begin{corollary}\label{cor:adjacency-permanent}
Under the assumptions of Proposition~\ref{prop:perfect-matching},
\begin{equation*}
\per(A(G))
=\per(B(G))\per(B(G)^T)
=\per(B(G))^2
=\mu(G)^2.
\end{equation*}
Consequently,
\begin{equation*}
\mu(G)=\sqrt{\per(A(G))}.
\end{equation*}
\end{corollary}

\subsection{Fibonacci, Lucas and Pairwise Dyadic Sequences}
\label{subsec:sequences}

The Fibonacci and Lucas sequences are well-known integer sequences with
important connections to number theory and other areas of mathematics
\cite{vajda2008}. Fibonacci numbers occur in biological growth patterns and plant phyllotaxis \cite{mitchison1977}, whereas Lucas sequences have applications in primality testing and cryptography \cite{ballot2023}. Koshy also discusses identities
that express Fibonacci and Lucas numbers through suitable sums of binomial
coefficients in Pascal's triangle \cite{koshy2011}.

\subsubsection*{The Fibonacci Sequence}
The Fibonacci sequence $(\mathfrak{f}_n)$ is defined by $\mathfrak{f}_n = \mathfrak{f}_{n-1}+\mathfrak{f}_{n-2}, \forall n\geq 2$, 
with initial conditions $\mathfrak{f}_0=0$ and $\mathfrak{f}_1=1$. Its initial terms are
\begin{equation*}
0,1,1,2,3,5,8,13,21,34,55,\ldots
\end{equation*}
and this sequence appears in the OEIS as A000045 \cite{OEIS}.  The following lemma\footnote{A standard identity that can be found in \cite{vajda2008}.} will be useful in the next section.

\begin{lemma}\label{Fibonacci sum}
For $n\ge 2$,  we have $\displaystyle\sum_{\ell=1}^n \mathfrak{f}_\ell=\mathfrak{f}_{n+2}-1$.
\end{lemma}

\begin{proof}
Proof is by induction on $n\ge 2$.  When $n=2$, we have $\mathfrak{f}_1+\mathfrak{f}_2=1+1=2=\mathfrak{f}_4-2$.  Assuming that the result holds for arbitrary $n\ge 2$, consider the sum over $\ell$ from $1$ to $n+1$:\begin{align*}
\sum_{\ell=1}^{n+1}\mathfrak{f}_\ell&=\sum_{\ell=1}^n\mathfrak{f}_\ell+\mathfrak{f}_{n+1}\\
&=(\mathfrak{f}_{n+2}-1)+\mathfrak{f}_{n+1}=\mathfrak{f}_{n+3}-1.
\end{align*}
\end{proof}

\subsubsection*{The Lucas Sequence}
The Lucas sequence $(\mathfrak{u}_n)$ is defined by $\mathfrak{u}_n=\mathfrak{u}_{n-1}+\mathfrak{u}_{n-2}, \forall n\geq 2$,
with initial conditions $\mathfrak{u}_0=2$ and $\mathfrak{u}_1=1$. Its initial terms are
\begin{equation*}
2,1,3,4,7,11,18,29,47,76,\ldots
\end{equation*}
and this sequence appears in the OEIS as A000032 \cite{OEIS}.  

\begin{lemma}\label{Fibonacci Lucas}
In terms of the Fibonacci sequence, the Lucas numbers satisfy $\mathfrak{u}_n=\mathfrak{f}_{n+1}+\mathfrak{f}_{n-1}$ 
for $n\ge 1$.
\end{lemma}

\begin{proof}
Again proof is by induction on $n\ge 1$.  When $n=1$, we have $\mathfrak{u}_1=1=1+0=\mathfrak{f}_2+\mathfrak{f}_0$.  Assuming the result holds for $n\ge 1$, we have \begin{align*}
\mathfrak{u}_{n+1}=\mathfrak{u}_{n}+\mathfrak{u}_{n-1}&=\mathfrak{f}_{n+1}+\mathfrak{f}_{n-1}+\mathfrak{f}_{n}+\mathfrak{f}_{n-2}\\
&=(\mathfrak{f}_{n+1}+\mathfrak{f}_{n})+(\mathfrak{f}_{n-1}+\mathfrak{f}_{n-2})\\
&=\mathfrak{f}_{n+2}+\mathfrak{f}_n.
\end{align*}
\end{proof}

\subsubsection*{Other Sequences of Interest}

In addition to Fibonacci and Lucas sequences, we consider the following sequences.

\begin{itemize} 
\item The sequence $\left(2^{\lfloor n/2\rfloor}\right)$.  Its terms are $(1,2,2,4,4,8,8,16,16,32,32,\ldots)$, and it can be found in the OEIS as sequence A060546 \cite{OEIS}.  

\item The sequence of Leonardo numbers, $(\mathfrak{d}_n)$, is defined by $\mathfrak{d}_1=1$, $\mathfrak{d}_2=3$, and $\mathfrak{d}_n=1+\mathfrak{d}_{n-1}+\mathfrak{d}_{n-2}$ for $n\ge 3$.  Leonardo numbers appear as sequence A001595 at \cite{OEIS}.

\item The sequence $(\mathfrak{a}_n)$ is the integer sequence defined by $\mathfrak{a}_0=1$, $\mathfrak{a}_1=1$, $\mathfrak{a}_2=3$, $\mathfrak{a}_3=4$, and $\mathfrak{a}_n=1+\mathfrak{a}_{n-1}+\mathfrak{a}_{n-2}$ for $n>3$.    Its initial terms are $ (1, 3, 4, 8, 13, 22, 36, 59, 96, 156, 253, 410, 664, 1075, 1740, 2816,\ldots)$.  

\end{itemize}

\section{Bipartite Matchings and Biadjacency Matrices}\label{matchings and matrices}

This section presents four results concerning the permanents of particular
families of biadjacency matrices and the numbers of perfect matchings in
their corresponding bipartite graphs. These results are established using
the permanent expansions and recurrence methods described in the preceding
section.
\begin{theorem}\label{thm:Cn}
Let $\mathcal{C}_n$ be the bipartite graph whose $n\times n$ biadjacency matrix is
$C_n=(c_{ij})$, where
\begin{equation*}
C_n=
\begin{pmatrix}
1 & 1 & 1 & \cdots & 1\\
1 & 1 & 0 & \cdots & 0\\
1 & 0 & 1 & \ddots & \vdots\\
\vdots & \vdots & \ddots & \ddots & 0\\
1 & 0 & \cdots & 0 & 1
\end{pmatrix},
\end{equation*}
and
\begin{equation*}
c_{ij}=
\begin{cases}
1, & \text{if }i=1,\\
1, & \text{if }i\geq2\text{ and }(j=1\text{ or }j=i),\\
0, & \text{otherwise}.
\end{cases}
\end{equation*}
Then, for every $n\geq1$, the number of perfect matchings of $\mathcal{C}_n$ is
$n$.
\end{theorem}

\begin{proof}
We proceed by induction on $n$. For $n=1$, we have $C_1=(1)$, and
hence $\per(C_1)=1$. Assume that
$\per(C_n)=n$ for some $n\geq1$.

Consider $C_{n+1}$. Its last column contains exactly two nonzero
entries, namely $c_{1,n+1}=1$ and $c_{n+1,n+1}=1$. Expanding the
permanent along the last column gives
\begin{equation*}
\per(C_{n+1})
=\per(M_1)+\per(C_n),
\end{equation*}
where $M_1$ is obtained by deleting the first row and the last column
of $C_{n+1}$. The second minor, obtained by deleting the last row and
the last column, is precisely $C_n$.

The last row of $M_1$ has only one nonzero entry, located in its first
column. Expanding along this row leaves the identity matrix
$I_{n-1}$; therefore, $\per(M_1)=1$. By the induction
hypothesis,
\begin{equation*}
\per(C_{n+1})
=1+\per(C_n)=1+n=n+1.
\end{equation*}
Thus, $\per(C_n)=n$ for every $n\geq1$. By
Proposition~\ref{prop:perfect-matching}, the number of perfect
matchings of $\mathcal{C}_n$ is $n$.
\end{proof}

\begin{figure}[h]
$$\includegraphics[width=240pt]{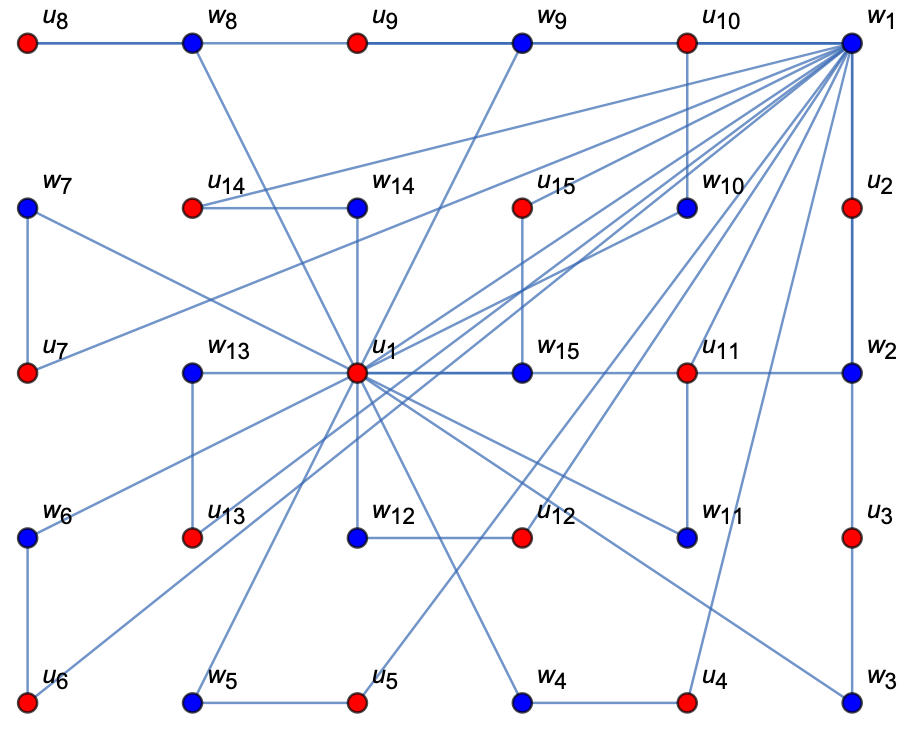}$$
\caption{The bipartite graph $\mathcal{C}_{15}$ corresponding to the
biadjacency matrix $C_{15}$ of Theorem \ref{thm:Cn}.}
\label{fig:C15}
\end{figure}

\begin{theorem}\label{thm:Wn}
Let $\mathcal{W}_n$ be the bipartite graph whose $n\times n$
biadjacency matrix is $W_n=(w_{ij})$, where
\begin{equation*}
w_{ij}=
\begin{cases}
1, & \text{if }j=1\text{ or }|i-j|\leq1,\\
0, & \text{otherwise}.
\end{cases}
\end{equation*}
Then, for every $n\geq2$, the number of perfect matchings of
$\mathcal{W}_n$ is $F_{n+2}-1$.
\end{theorem}

\begin{proof}
Let $W_n(i\mid1)$ denote the matrix obtained from $W_n$ by deleting
its $i$-th row and first column. Since every entry in the first
column of $W_n$ is $1$, expanding the permanent along that column
gives
\begin{equation*}
\per(W_n)
=\sum_{i=1}^{n}\per\bigl(W_n(i\mid1)\bigr).
\end{equation*}

For each $i$, the first $i-1$ rows of $W_n(i\mid1)$ force the
successive selection of their available diagonal entries. After
deleting the corresponding rows and columns, the remaining
nontrivial block is the $(n-i)\times(n-i)$ tridiagonal matrix whose
main, super-, and subdiagonal entries are all $1$. By Theorem~2.4 of
\"{O}tele\c{s}~\cite{oteles2017},
\begin{equation*}
\per\bigl(W_n(i\mid1)\bigr)=\mathfrak{f}_{n-i+1}.
\end{equation*}
For $i=n$, this uses the convention that the permanent of the empty
matrix is $1=\mathfrak{f}_1$.  In light of Lemma \ref{Fibonacci sum}, we have
\begin{align*}
\per(W_n)=\sum_{i=1}^{n}\mathfrak{f}_{n-i+1}&=\sum_{k=1}^{n}\mathfrak{f}_k\\
&=\mathfrak{f}_{n+2}-1.
\end{align*}
Therefore, by Proposition~\ref{prop:perfect-matching}, the number of
perfect matchings of $\mathcal{W}_n$ is $\mathfrak{f}_{n+2}-1$.
\end{proof}

\begin{figure}[h]
\centering
\includegraphics[width=240pt]{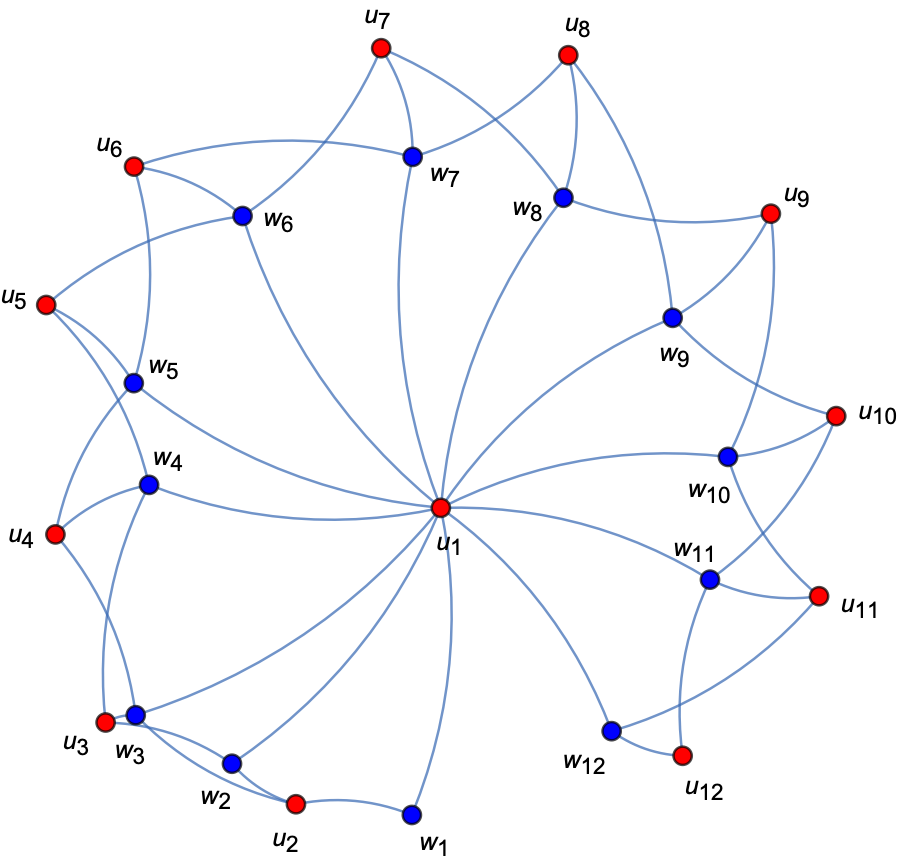}
\caption{The bipartite graph $\mathcal{W}_{12}$ corresponding to the
biadjacency matrix $W_{12}$ of Theorem \ref{thm:Wn}.}
\label{fig:W12}
\end{figure}

\begin{theorem}\label{thm:Rn}
For $n\in\mathbb{N}$, let $\mathcal{R}_n$ be the bipartite graph on $2n$ vertices $V=U\sqcup W$, where $|U|=|W|=n$, and whose biadjacency matrix $R_n$ is given by  
$R_1=(1)$ and  
\begin{equation*}
R_n=\begin{cases}
\begin{pmatrix}
1&1\\
1&1
\end{pmatrix}\otimes \mI_{n/2}&\text{\rm if }n\equiv 0\pmod{2},\\[1em]
\begin{pmatrix}
R_{n-1}&0\\
0& 1\end{pmatrix}&\text{\rm if }n\equiv 1\pmod{2},
\end{cases}
\end{equation*}
for $n>1$.  Here, $\mI_{n/2}$ denotes the identity matrix of order $n/2$. Then, for every $n\geq1$, the number of perfect matchings of $\mathcal{R}_n$ is
\begin{equation*}
\per(R_n)
=2^{\lfloor n/2\rfloor}
\end{equation*}
where $(2^{\lfloor n/2\rfloor})_{n\geq1}$ is the pairwise dyadic sequence.
\end{theorem}

\begin{figure}[h]
$$\includegraphics[width=230pt]{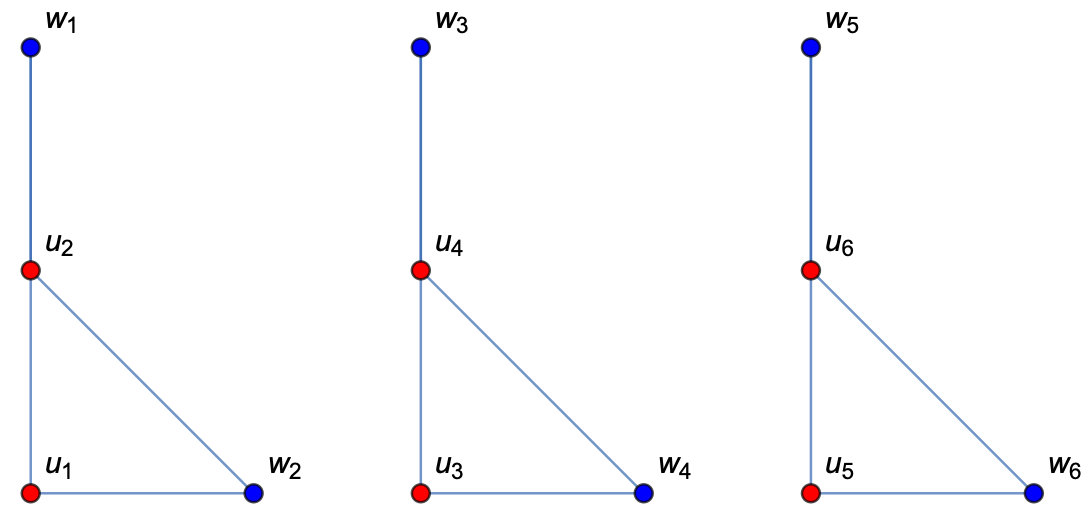}$$
\caption{The bipartite graph $\mathcal{R}_6$ corresponding to the
biadjacency matrix $R_6$.}
\label{fig:R6}
\end{figure}

\begin{proof}
Since $\per(J)=2$ and $\per((1))=1$, we compute the permanent separately for even and odd values of $n$.  The permanent of a block-diagonal matrix is the product of the
permanents of its diagonal blocks because every nonzero term in its permanent expansion selects entries entirely within those blocks.

If $n=2m$, then
\begin{equation*}
\per(R_{2m})
=\prod_{k=1}^{m}\per(J)
=2^m.
\end{equation*}
If $n=2m+1$, then
\begin{equation*}
\per(R_{2m+1})
=\left(\prod_{k=1}^{m}\per(J)\right)
\per((1))
=2^m.
\end{equation*}
Therefore, in both cases,
\begin{equation*}
\per(R_n)=2^{\lfloor n/2\rfloor}=D_n.
\end{equation*}
By Proposition~\ref{prop:perfect-matching}, the number of perfect
matchings of $\mathcal{R}_n$ is $D_n$.
\end{proof}

\begin{theorem}\label{thm:Qn}
Let $\mathcal{Q}_n$ be the bipartite graph whose $n\times n$
biadjacency matrix is $Q_n=(q_{ij})$, where
\begin{equation*}
q_{ij}=\begin{cases}
1, & i=1\text{ and }
j\in\{1,2\}\cup\{4,5,\ldots,n\},\\
1&2\le i<n\,\text{\rm and }j\in\{i-1,i,i+1\},\\
1&i=n \,\text{\rm and }j\in\{n-1,n\},\\
0&\text{\rm otherwise.}
\end{cases}
\end{equation*}
Then, for every $n\geq2$, the number of perfect matchings of
$\mathcal{Q}_n$ is $\mathfrak{u}_n-1$.
\end{theorem}

\begin{proof}
Direct calculation gives
\begin{equation*}
\per(Q_2)=2=\mathfrak{u}_2-1
\quad\text{and}\quad
\per(Q_3)=3=\mathfrak{u}_3-1.
\end{equation*}
Now let $n\geq4$. For $1\leq j\leq n$, let
$Q_n(1\mid j)$ denote the submatrix obtained by deleting the first
row and the $j$-th column of $Q_n$. Expanding the permanent along
the first row gives
\begin{equation*}
\per(Q_n)
=
\per\bigl(Q_n(1\mid1)\bigr)
+\per\bigl(Q_n(1\mid2)\bigr)
+\sum_{j=4}^{n}\per\bigl(Q_n(1\mid j)\bigr).
\end{equation*}

By Theorem 2.4 of \"{O}tele\c{s} \cite{oteles2017}, together with invariance
of the permanent under row and column permutations,
\begin{equation*}
\per\bigl(Q_n(1\mid1)\bigr)=\mathfrak{f}_n
\end{equation*}
and
\begin{equation*}
\per\bigl(Q_n(1\mid2)\bigr)=\mathfrak{f}_{n-1}.
\end{equation*}
For $4\leq j\leq n$, successive expansions along rows or columns
containing a unique nonzero entry reduce $Q_n(1\mid j)$ to the
corresponding matrix of order $n-j$. Therefore,
\begin{equation*}
\per\bigl(Q_n(1\mid j)\bigr)=\mathfrak{f}_{n-j+1}.
\end{equation*}
Consequently,
\begin{align*}
\per(Q_n)
&=\mathfrak{f}_n+\mathfrak{f}_{n-1}
  +\sum_{j=4}^{n}\mathfrak{f}_{n-j+1}\\
&=\mathfrak{f}_n+\mathfrak{f}_{n-1}
  +\sum_{k=1}^{n-3}\mathfrak{f}_k.
\end{align*}
Using Lemmas \ref{Fibonacci sum} and \ref{Fibonacci Lucas}, we obtain
\begin{align*}
\per(Q_n)
&=\mathfrak{f}_n+2\mathfrak{f}_{n-1}-1\\
&=\mathfrak{f}_{n+1}+\mathfrak{f}_{n-1}-1\\
&=\mathfrak{u}_n-1.
\end{align*}
Therefore, by Proposition~\ref{prop:perfect-matching}, the number of
perfect matchings of $\mathcal{Q}_n$ is $\mathfrak{u}_n-1$.
\end{proof}

\begin{figure}[ht]
\centering
\includegraphics[width=240pt]{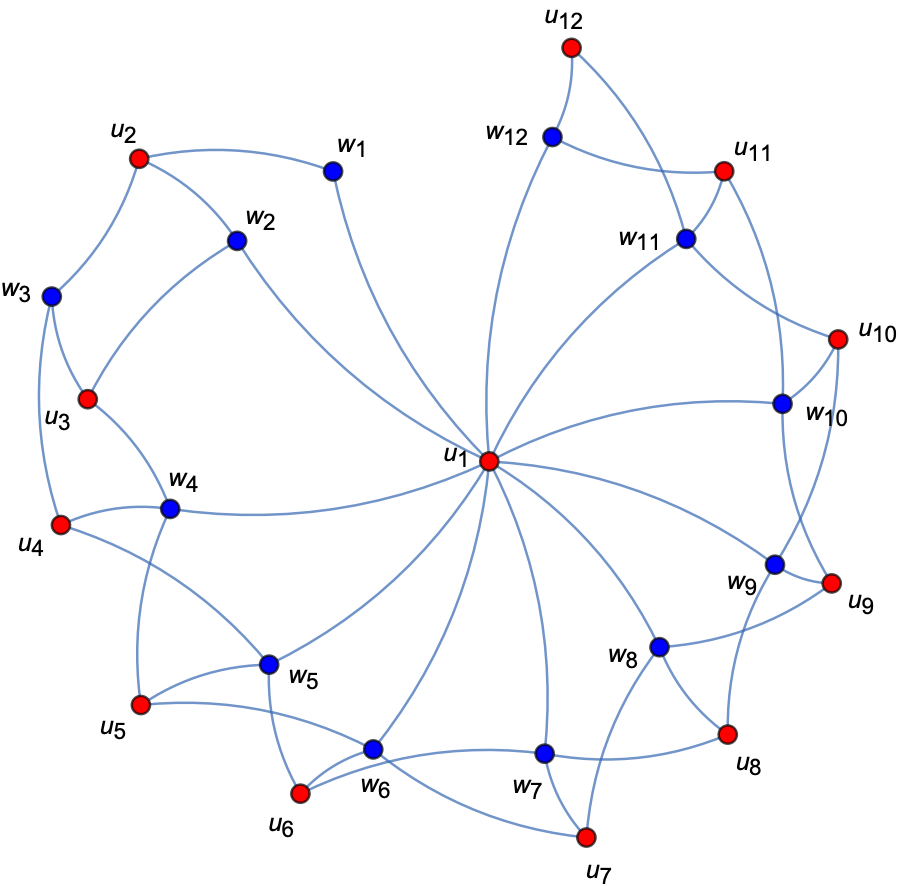}
\caption{The bipartite graph $\mathcal{Q}_{12}$ corresponding to the
biadjacency matrix $Q_{12}$.}
\label{fig:Q10}
\end{figure}

\begin{theorem}\label{Leonardo}
Let $\mathcal{B}_n$ be the bipartite multigraph on $2n$ vertices whose biadjacency matrix $B_n=(b_{ij})$ is defined by 
\begin{equation*}
b_{ij}=
\begin{cases}
1, & \text{if }i=1\text{ and }
j\in\{1\}\cup\{3,4,\ldots,n\},\\
2, & \text{if }i=1, j=2\\
1, & \text{if }2\leq i\leq n\text{ and }
\vert i-j\vert\le 1,\\
0, & \text{otherwise}.
\end{cases}
\end{equation*}
Then, for every $n\geq2$, the number of perfect matchings of
$\mathcal{B}_n$ is the $n$th Leonardo number $\mathfrak{d}_n$, where $(\mathfrak{d}_n)$ is defined by $\mathfrak{d}_1=1$, $\mathfrak{a}_2=3$, $\mathfrak{a}_3=5$, $\mathfrak{a}_n=1+\mathfrak{a}_{n-1}+\mathfrak{a}_{n-2}$ for $n\ge 4$.
\end{theorem}

\begin{proof}
The first four matrices are \begin{align*}
B_1&=\begin{pmatrix} 1 \end{pmatrix}& B_2&=\begin{pmatrix}1&2\\1&1\end{pmatrix}\\
B_3&=\begin{pmatrix}1&2&1\\1&1&1\\0&1&1\end{pmatrix}&B_4&=\begin{pmatrix}1&2&1&1\\1&1&1&0\\0&1&1&1\\0&0&1&1\end{pmatrix}.
\end{align*}
By construction, when $n\ge 4$, the $n$th matrix is of the form \begin{equation*}
B_n\hskip10pt=\hskip10pt\begin{pmatrix}
&&&&1\\
&&&&0\\
&B_{n-1}&&&0\\
&&&&0\\
&&&&1\\
0&\cdots&0&1&1
\end{pmatrix}\hskip10pt=\hskip10pt\begin{pmatrix}
&&&1&1\\
B_{n-2}&&&0&0\\
&&&0&\vdots\\
&&&1&0\\
\cdots&0&1&1&1\\
\cdots&0&0&1&1
\end{pmatrix}.
\end{equation*}
Clearly, $\per(B_1)=1$, $\per(B_2)=3$, and $\per(B_3)=5$, establishing the first three terms of the sequence.  When $n>3$, we compute $\per(B_n)$ by Laplace expansion in the last row, giving \begin{align*}
\per(B_n)&=\per(B_n(n\mid n-1))+\per(B_n(n\mid n))\\
&=\per\begin{pmatrix}
&&&1\\
B_{n-2}&&&0\\
&&&\vdots\\
&&&0\\
\cdots&0&1&1
\end{pmatrix}+\per(B_{n-1}).
\end{align*}
Computing the remaining permanent by Laplace expansion in the last column results in the sum of two permanents: the first is that of an upper triangular matrix with 1s down the main diagonal and the second is that of $B_{n-2}$.  Hence,
\begin{align*}
\per\begin{pmatrix}
&&&1\\
B_{n-2}&&&0\\
&&&\vdots\\
&&&0\\
\cdots&0&1&1
\end{pmatrix}&=1+\per(B_{n-2}).
\end{align*}  
Finally, we arrive at our result: $\per(B_n)=1+\per(B_{n-1})+\per(B_{n-2})$.  The resulting sequence of determinants is then identical to the sequence $(\mathfrak{d}_n)$.
\end{proof}

The next family of bipartite multigraphs results in a slightly different integer sequence satisfying the same recurrence relation.

\begin{theorem}\label{thm:Kn}
Let $\mathcal{K}_n$ be the bipartite multigraph on $2n$ vertices whose biadjacency matrix $K_n=(k_{ij})$ is defined by 
\begin{equation*}
k_{ij}=
\begin{cases}
1, & \text{if }i=1\text{ and }
j\in\{1\}\cup\{4,\ldots,n\},\\
2, & \text{if }i=1, j=2\\
1, & \text{if }2\leq i\leq n\text{ and }
\vert i-j\vert\le 1,\\
0, & \text{otherwise}.
\end{cases}
\end{equation*}
Then, for every $n\geq2$, the number of perfect matchings of
$\mathcal{K}_n$ is $\mathfrak{a}_n$, where $(\mathfrak{a}_n)$ is defined by $\mathfrak{a}_1=1$, $\mathfrak{a}_2=3$, $\mathfrak{a}_3=4$, $\mathfrak{a}_n=1+\mathfrak{a}_{n-1}+\mathfrak{a}_{n-2}$ for $n\ge 4$.
\end{theorem}

\begin{proof}
The first four matrices are \begin{align*}
K_1&=\begin{pmatrix} 1 \end{pmatrix}& K_2&=\begin{pmatrix}1&2\\1&1\end{pmatrix}\\
K_3&=\begin{pmatrix}1&2&0\\1&1&1\\0&1&1\end{pmatrix}&K_4&=\begin{pmatrix}1&2&0&1\\1&1&1&0\\0&1&1&1\\0&0&1&1\end{pmatrix}.
\end{align*}
By construction, when $n\ge 4$, the $n$th matrix is of the form \begin{equation*}
K_n\hskip10pt=\hskip10pt\begin{pmatrix}
&&&&1\\
&&&&0\\
&K_{n-1}&&&0\\
&&&&0\\
&&&&1\\
0&\cdots&0&1&1
\end{pmatrix}\hskip10pt=\hskip10pt\begin{pmatrix}
&&&1&1\\
K_{n-2}&&&0&0\\
&&&0&\vdots\\
&&&1&0\\
\cdots&0&1&1&1\\
\cdots&0&0&1&1
\end{pmatrix}.
\end{equation*}
Clearly, $\per(K_1)=1$, $\per(K_2)=3$, and $\per(K_3)=4$, establishing the first three terms of the sequence.  When $n>3$, we compute $\per(K_n)$ by Laplace expansion in the last row, as in the proof of Theorem \ref{Leonardo}, giving $\per(K_n)=1+\per(K_{n-1})+\per(K_{n-2})$.  Letting $\mathfrak{a}_n=\per(K_n)$ completes the proof.
\end{proof}

\section{Zeon Algebra and Bipartite Matchings}\label{zeon formulations}

Zeon algebras can be interpreted as commutative fermion algebras.  Equivalently, they can be constructed as subalgebras of Clifford algebras.  They provide a natural language for problems in combinatorics and algebraic probability~\cite{CAZ}.  Zeons and permanents have been discussed in work by Feinsilver and McSorley \cite{FeinsilverMcSorleyIJC}.  Some connections between integer sequences and zeons have also been explored in work by Neto \cite{NetoJIS2014,NetoJIS2015}.

\subsection{Zeon Preliminaries}
For $n\in\mathbb{N}$,  let $\CZ_n$ denote the complex commutative algebra generated by the collection $\{\zeta_{\{i\}}: 1\le i\le n\}$ along with the scalar $1=\zeta_\varnothing$ subject to the following multiplication rules:
\begin{gather*}
\zeta_{\{i\}}\,\zeta_{\{j\}}=\zeta_{\{i,j\}}=\zeta_{\{j\}}\,\zeta_{\{i\}}\,\,\text{ \rm
for }i\ne j,\text{ \rm and}\\
{\zeta_{\{i\}}}^2=0\,\,\text{ \rm for }1\le i\le n.
\end{gather*}
It is evident that a general element $u\in\CZ_n$ can be expanded as $u=\displaystyle\sum_{I\in2^{[n]}}u_{I}\,\zeta_{I}$, or more simply as $\sum_I u_I\zeta_I$,  where $I\in2^{[n]}$ is a subset of the {\em $n$-set}, $[n]:=\{1,2,\ldots,n\}$, used as a multi-index, $u_{I}\in\mathbb{R}$, and $\displaystyle\zeta_{I}=\prod_{\iota\in I}\zeta_\iota$.  The algebra $\CZ_n$ is called the ($n$-particle) {\em complex zeon algebra}.  
As a vector space, this $2^n$-dimensional algebra has a canonical basis of  {\em basis blades} of the form $\{\zeta_I: I\subseteq [n]\}$.  The null-square property of the generators $\{\zeta_i:1\le i\le n\}$ guarantees that the product of two basis blades satisfies the following:\begin{equation}\label{blade product}
\zeta_I\zeta_J=\begin{cases}
\zeta_{I\cup J}& I\cap J=\varnothing,\\
0&\text{\rm otherwise.}
\end{cases}
\end{equation}
It should be clear that $\CZ_n$ is graded.  For non-negative integer $k$, the {\em grade-$k$ part} of element $u=\sum_I u_I\zeta_I$ is defined as \begin{equation*}
\langle u\rangle_k=\sum_{\{I:|I|=k\}}u_I\zeta_I.
\end{equation*}
An inner product is defined on $\CZ_n$ by linear extension of 
\begin{equation*}
\displaystyle\left<\sum_I u_I\zeta_I,\zeta_J\right>=u_J.
\end{equation*}

The basis blades $\{\zeta_I: I\in 2^{[n]}\}$ form an orthonormal collection with respect to this inner product; thus, any element $u\in\CZ_n$ can be expanded as $\displaystyle u=\sum_{I}\left<u,\zeta_I\right>\zeta_I.$   Moreover, an inner product  norm is naturally induced on $\CZ_n$ by $\|u\|=\langle u,u\rangle^{1/2}$.    

It is often convenient to separate the scalar term of a zeon element from the rest of it. To this end, for $z\in\CZ_n$ we write $\zco z=\langle z\rangle_0$ for the {\em scalar part} of $z$, and $\mathfrak{D}z=z-\zco z$ for the {\em dual part} of $z$.

The finite-dimensional algebra is extended to zeon algebras of arbitrary dimension by defining \begin{equation*}
\CZ={\rm span}\{\zeta_I: I\subset\mathbb{N}, I \,\text{\rm finite}\}, 
\end{equation*}
which we call the {\em infinite-dimensional (complex) zeon algebra}.

\subsubsection{Zeon Formulations of Bipartite Matching Problems}

\begin{definition}
Let $G=(V,E)$ be a simple graph on $n$ labeled vertices $V=\{v_1, \ldots, v_n\}$ and edges $E$ consisting of unordered pairs of vertices.  The {\em nilpotent adjacency matrix of $G$} is defined as the $n\times n$ matrix $\Phi=(\varphi_{ij})$ having entries in ${\CZ_n}$ determined as follows:
\begin{equation*}
\varphi_{ij}=\begin{cases}
0&\{v_i, v_j\}\notin E,\\
\zeta_{\{j\}}& \{v_i, v_j\}\in E.
\end{cases}
\end{equation*}
Equivalently, $\Phi=A\mathcal{Z}_n$, where $A$ is the graph's adjacency matrix and $\mathcal{Z}_n={\rm diag}(\zeta_{\{1\}}, \ldots, \zeta_{\{n\}})$ is the diagonal matrix of zeon generators.
\end{definition}

Modifying Dirac notation, let us define $\langle \zeta_{\{i\}} \vert$ to denote the row vector whose entries are all zero, except the $i$th column, which is $\zeta_{\{i\}}$.  Similarly, $\vert j\rangle$ will denote the column vector with $1$ in the $j$th row and zeros elsewhere.  

Because the generators of $\CZ_n$ square to zero, a straightforward inductive argument shows that the nonzero terms of $\left<v_i\vert\Phi^k\vert v_j\right>$ are multivectors corresponding to two types of $k$-walks from $v_i$ to $v_j$:
self-avoiding walks (i.e., walks with no repeated vertices) and walks in which $v_i$ is repeated exactly once at some step but are otherwise self-avoiding.   For $i\ne j$, multiplication by $\zeta_{\{i\}}$ cancels walks of the second type, so that nonzero terms represent paths from $v_i$ to $v_j$.   

When $i=j$, walks of the second type are zeroed in the $k^{\rm th}$ step when the walk is closed. Hence, terms of $\left<v_i\vert \Phi^k\vert v_i\right>$ represent the collection of $k$-cycles based at $v_i$.  We thereby arrive at the following result as formulated in \cite{zeon laplacian}.

\begin{description}
\item [Nil-Structure Theorem.]
Let $\Phi$ be the nilpotent adjacency matrix of an
$n$-vertex graph $G$, and let $k\in\mathbb{N}$.  For $1\le i\ne j \le n$, 
\begin{equation*}
\left<\zeta_{\{i\}}\vert\Phi^k\vert v_j\right>=\sum_{\genfrac{}{}{0pt}{}{(w_1,\ldots, w_k)\in V^k}{(w_k=v_j)\wedge(m\ne \ell\Rightarrow w_m\ne w_\ell)}}\zeta_{(w_1,\ldots, w_k)}=
\sum_{\genfrac{}{}{0pt}{}{I\subseteq V}{|I|=k+1}}\omega_I \zeta_{I},\end{equation*}
where $\omega_I$ denotes the number of $k$-paths from $v_i$ to $v_j$.  Further,  when $1\le i \le n$,
\begin{equation*}
\left<v_i\vert\Phi^k\vert v_i\right>=\sum_{\genfrac{}{}{0pt}{}{I\subseteq V}{|I|=k}}\omega_I \zeta_{I},\end{equation*}
where $\omega_I$ denotes the number of $k$-cycles on vertex set $I$ based at $v_i\in I$.  
\end{description}

The zeon bipartite matchings theorem can now be presented.

\begin{theorem}[Zeon Bipartite Matchings Theorem]\label{zeon matchings theorem}
Let $A$ be the adjacency matrix of a graph $G=(V,E)$ on $|V|=2n$ vertices, and let $\Phi=A\mathcal{Z}_{2n}$ be the nilpotent adjacency matrix of $G$.  If $\tr(\zexp(\Phi))\in\langle \CZ\rangle_{\rm even}$, then $G$ is bipartite.  Moreover, letting $\gamma=\frac{1}{2}\tr (\Phi^2)$, we have 
\begin{equation*}
\langle\zexp(\gamma), \zeta_{[2n]}\rangle=\left({\rm per}\,A\right)^{1/2}=\mu_G,
\end{equation*}
where $\mu_G$ is the number of perfect matchings of $G$.  
\end{theorem}

\begin{proof}
By the nil-structure theorem, $\gamma=\frac{1}{2}\tr (\Phi^2)$ is a sum of grade-2 basis blades representing the edges (2-cycles) of the graph.  Since $\zeta_{\{i,j\}}\zeta_{\{k,\ell\}}\ne 0$ if and only if $\{i,j\}\cap\{k,\ell\}=\varnothing$, we see that   
\begin{equation}
\gamma^n=n!\,\mu_G\,\zeta_{[2n]},
\end{equation} 
where $\mu_G$ is the number of perfect matching of the graph $G$.  
\end{proof}

More generally, the number of perfect matchings in any graph $G$ is given by powers of $\gamma=\frac{1}{2}\tr(\Phi^2)$ by removing the requirement that the exponential be in the even subalgebra.  Hence, we obtain the following corollary.

\begin{corollary}[Zeon Matchings]\label{zeon matchings corollary}
Let $A$ be the adjacency matrix of a graph $G=(V,E)$ on $|V|=2n$ vertices, and let $\Phi=A\mathcal{Z}_{2n}$ be the nilpotent adjacency matrix of $G$.  Letting $\gamma=\frac{1}{2}\tr (\Phi^2)$, we have 
\begin{equation*}
\langle\zexp(\gamma), \zeta_{[2n]}\rangle=\mu_G,
\end{equation*}
where $\mu_G$ is the number of perfect matchings of $G$.  
\end{corollary}

Letting $\Phi$ be the nilpotent adjacency matrix of a simple graph $G$ on $n$ vertices, the unique zeon element $\gamma=\tr(\Phi^2)/2$ will be referred to as the {\em zeon matching representation} of $G$.  Noting that $\gamma$ is a sum over the edges of $G$, it is a simple matter to reconstruct the adjacency matrix $A=(a_{ij})$ of $G$ by setting 
$a_{ij}= \langle \gamma, \zeta_{\{i,j\}}\rangle$ for each $1\le i\ne j\le n$ and $a_{jj}=0$ for $j=1, \ldots, n$.

The following lemma will be particularly useful for the next section, where we recover results using formal power series and zeon sequences.

\begin{lemma}[$\gamma$ from $B$]\label{biadjacency zeon lemma}
Let $B=(b_{ij})$ be the biadjacency matrix of a balanced bipartite graph $G$ on $2n$ vertices.  A zeon matching representation of $G$ is then given by \begin{equation}\label{biadjacency matching rep}
\gamma=\sum_{1\le i,j\le n} b_{ij}\zeta_{\{2i-1,2j\}}.
\end{equation}
In particular, if $\mu$ is the number of perfect matchings of $G$, then \begin{equation*}
\mu=\per(B)=\frac{1}{n!}\langle \gamma^{n}, \zeta_{[2n]}\rangle.
\end{equation*}
\end{lemma}

\begin{proof}
In light of previous results, it is sufficient to prove that $\gamma$ is a sum over the edges of $G$.  Given the bipartition $V=U\sqcup W$, we let vertices in $U$ be indexed by odd integers $\{1, 3, \ldots, 2n-1\}$ and vertices of $W$ be indexed by $\{2, 4, \ldots, 2n\}$.  It is then clear that the sum appearing in  \eqref{biadjacency matching rep} is a sum over the edges of $G$.
\end{proof}

\section{Zeon Sequences and the $\Omega$-Exponential}\label{OmegaExp}

Given a zeon sequence $\Gamma=(\gamma_n)$, we define the formal power series 
\begin{equation*}
\Wexp(t\,\Gamma)=\sum_{n=0}^\infty \frac{t^n{\gamma_n}^n}{(n!)^2},
\end{equation*}
and refer to this as the {\em $\Omega$-exponential of the zeon sequence $t\,\Gamma=(t\gamma_n)$}.

As a formal power series in the scalar variable $t$, the $k$th derivative of $\Wexp(t\,\Gamma)$ is given by 
\begin{align*}
\frac{d^k}{dt^k}\Wexp(t\Gamma)=k!\sum_{n=k}^\infty\frac{t^{n-k}{\gamma_n}^n}{(n!)^2}
&=\sum_{n=k}^\infty\frac{t^{n-k}{\gamma_n}^n}{n!(n-k)!}\\
&=\sum_{n=0}^\infty\frac{t^{n}{\gamma_{n+k}}^{n+k}}{n!(n+k)!}\\
&=\frac{{\gamma_k}^{k}}{k!}+\sum_{n=1}^\infty\frac{t^{n}{\gamma_{n+k}}^{n+k}}{n!(n+k)!}.
\end{align*}
From this, the $n$th term of the sequence is formally extracted by considering a limit as $t$ approaches zero.

\begin{proposition}[Pairwise dyadic sequence\footnote{The pairwise dyadic sequence is  sequence A016116 in the OEIS \cite{OEIS}.}]\label{pairwise dyadic zeon sequence}
Let the zeon sequence $\Gamma_{\mathcal{R}}=(\rho_n)$ be defined recursively  by $\rho_0=0$, and \begin{equation}\label{dyadic zeon recurrence}
\rho_n=\begin{cases}
\rho_{n-2}+\zeta_{\{2n-2,2n-1\}}+\zeta_{\{2n-1,2n\}}+\zeta_{\{2n-3,2n-2\}}+\zeta_{\{2n-3,2n\}}& n \text{ \rm even},\\
\rho_{n-1}+\zeta_{\{2n-1,2n\}}&n\text{ \rm odd},
\end{cases}
\end{equation}  
for positive integers $n$. Then, \begin{equation*}
\left< \lim_{t\to 0}\frac{d^n}{dt^n}\Wexp(t\Gamma_{\mathcal{R}}), \zeta_{[2n]}\right>= 2^{\lfloor n/2\rfloor}.
\end{equation*}
\end{proposition}

\begin{proof}
By Lemma \ref{biadjacency zeon lemma}, $\rho_n$ is, by construction, a zeon matching representation of the graph whose biadjacency matrix is $R_n$, as seen in Theorem \ref{thm:Rn}.  Hence, the result.
\end{proof}

\begin{example}
Let $G=(V,E)$ be the bipartite graph on $2n$ vertices $V=U\sqcup W$, where $|U|=|W|=n$ and whose biadjacency matrix $R_n$ is given by  \begin{equation*}
R_n=\begin{cases}
\begin{pmatrix}
1&1\\
1&1
\end{pmatrix}\otimes \mI_{n/2}&\text{\rm if }n\equiv 0\pmod{2},\\[1em]
\begin{pmatrix}
R_{n-1}&0\\
0& 1\end{pmatrix}&\text{\rm if }n\equiv 1\pmod{2}.
\end{cases}
\end{equation*}
The number of perfect matchings of $G$ is then\begin{equation*}
\mu_G=\left< \lim_{t\to 0}\frac{d^n}{dt^n}\Wexp(t\Gamma_{\mathcal{R}}), \zeta_{[2n]}\right>=2^{\lfloor n/2\rfloor},
\end{equation*}
where $\Gamma_{\mathcal{R}}$ is the zeon sequence defined in \eqref{dyadic zeon recurrence}.
\end{example}

\begin{proposition}[Fibonacci numbers\footnote{Fibonacci numbers are sequence A000045 in the OEIS \cite{OEIS}.}]\label{Fibonacci zeon sequence}
Let the zeon sequence $\Gamma_{\mathfrak{f}}=(\phi_n)$ be defined recursively by \begin{align*}
\phi_1&= \zeta_{\{1,2\}},\\
\phi_2&=\zeta_{\{1,2\}}+\zeta_{\{1,4\}}+\zeta_{\{2,3\}}+\zeta_{\{3,4\}},\\\phi_n&=\varphi_{n-1}+\zeta _{\{1,2n\}}+\zeta _{\{2n-3,2n\}}+\zeta _{\{2n-2,2n-1\}}+\zeta _{\{2n-1,2n\}} &(n\ge 3).
\end{align*}  
Then, \begin{equation*}
\left< \lim_{t\to 0}\frac{d^n}{dt^n}\Wexp(t\Gamma_{\mathfrak{f}}), \zeta_{[2n]}\right>= \mathfrak{f}_{n+2}-1,
\end{equation*}
where $(\mathfrak{f}_n)$ is the Fibonacci sequence defined by $\mathfrak{f}_0=1$, $\mathfrak{f}_1=1$, and $\mathfrak{f}_n=\mathfrak{f}_{n-1}+\mathfrak{f}_{n-2}$ for $n\ge2$.  
\end{proposition}

\begin{proof}
By construction, $\phi_n$ is a zeon matching representation of the graph whose biadjacency matrix is $W_n$.  The result then follows from Theorem \ref{thm:Wn}
\end{proof}

\begin{example}
Let $G=(U,V,E)$ be the bipartite graph on $2n$ vertices $U\sqcup V$, where $|U|=|V|=n$ and whose biadjacency matrix $W_n(w_{ij})$ is determined by \begin{equation*}
w_{ij}=\begin{cases}
1&j=1 \text{ \rm or }|i-j|\le 1,\\
0&\text{\rm otherwise.}
\end{cases}
\end{equation*}
Letting $\Gamma_{\mathfrak{f}}=(\phi_n)$ be the zeon sequence defined in Proposition \ref{Fibonacci zeon sequence}, the number of perfect matchings of $G$ is \begin{equation*}
\mu_G=\left< \lim_{t\to 0}\frac{d^n}{dt^n}\Wexp(t\Gamma_{\mathfrak{f}}), \zeta_{[2n]}\right>=\mathfrak{f}_{n+2}-1.
\end{equation*}
\end{example}

\begin{proposition}[Lucas numbers\footnote{Lucas numbers can be found as sequence A000032 in the OEIS \cite{OEIS}.}]\label{Lucas zeon sequence}
Let the zeon sequence $\Gamma_{\mathcal{L}}=(\lambda_n)$ be defined recursively by \begin{align*}
\lambda_1&=\zeta_{\{1,2\}},\\ 
\lambda_2&=\zeta_{\{1,2\}}+\zeta_{\{1,4\}}+\zeta_{\{2,3\}}+\zeta_{\{3,4\}},\\
\lambda_3&=\lambda_2+\zeta_{\{5,6\}}+\zeta_{\{4,5\}}+\zeta_{\{3,6\}},\\
\lambda_n&=\lambda_{n-1}+\zeta_{\{2n-1,2n\}}+\zeta_{\{1,2n\}}+\zeta_{\{2n-1,2n-2\}}+\zeta_{\{2n-3,2n\}},\,\,\, (n\ge 4).
\end{align*} 
Then, for every $n\ge 2$, \begin{equation*}
\left< \lim_{t\to 0}\frac{d^n}{dt^n}\Wexp(t\Gamma_{\mathcal{L}}), \zeta_{[2n]}\right>= \mathfrak{u}_n-1, 
\end{equation*}
where $(\mathfrak{u}_n)$ is the Lucas sequence defined by $\mathfrak{u}_0=2$, $\mathfrak{u}_1=1$, and $\mathfrak{u}_n=\mathfrak{u}_{n-1}+\mathfrak{u}_{n-2}$ for $n\ge 2$.
\end{proposition}

\begin{proof}
By construction, the sequence $\lambda_n$ is a zeon matching representation of the bipartite graph whose biadjacency matrix is $W_n$ seen in Theorem \ref{thm:Wn}.  Hence, the result.
\end{proof}

\begin{example}
Let $G=(U,V,E)$ be the bipartite graph on $2n$ vertices $U\sqcup V$, where $|U|=|V|=n$ and whose biadjacency matrix $Q_n=(q_{ij})$ is determined by \begin{equation*}
q_{ij}=\begin{cases}
1, & i=1\text{ and }
j\in\{1,2\}\cup\{4,5,\ldots,n\},\\
1&2\le i<n\,\text{\rm and }j\in\{i-1,i,i+1\},\\
1&i=n \,\text{\rm and }j\in\{n-1,n\},\\
0&\text{\rm otherwise.}
\end{cases}
\end{equation*}
Letting $\Gamma_{\mathcal{L}}=(\lambda_n)$ be the zeon sequence defined in Proposition \ref{Lucas zeon sequence}, the number of perfect matchings of $G$ is \begin{equation*}
\mu_G=\left< \lim_{t\to 0}\frac{d^n}{dt^n}\Wexp(t\Gamma_{\mathcal{L}}), \zeta_{[2n]}\right>=\mathfrak{u}_n-1
\end{equation*}
for every $n\ge 2$.
\end{example}

\begin{proposition}[Leonardo numbers\footnote{Leonardo numbers appear as sequence A001595 in the OEIS \cite{OEIS}.}]\label{Leonardo zeon sequence}
Let the zeon sequence $\Gamma_{\mathfrak{D}}=(\delta_n)$ be defined recursively by \begin{align}\label{zeon Leonardo}
\delta_1&=\zeta_{\{1,2\}}\nonumber\\ 
\delta_n&=\delta_{n-1}+\zeta_{\{2n-1,2n\}}+\zeta_{\{1,2n\}}+\zeta_{\{2n-1,2n-2\}}+\zeta_{\{2n-3,2n\}},\,\,\, (n\ge 2).
\end{align} 
Then,  \begin{equation*}
 \Wexp(t\,\Gamma_{\mathfrak{D}})=\sum_{n=0}^\infty t^n\mathfrak{d}_n\zeta_{[2n]}, 
\end{equation*}
where $\mathfrak{d}_n$ is the sequence of Leonardo numbers defined by $\mathfrak{d}_1=1$, $\mathfrak{d}_2=3$, and $\mathfrak{d}_n=1+\mathfrak{d}_{n-1}+\mathfrak{d}_{n-2}$ for $n\ge 3$.
In particular, for positive integer $k$ we have
\begin{equation*}
\lim_{t\to 0}\frac{d^k}{dt^k}\Wexp(t\,\Gamma_{\mathfrak{D}})=\mathfrak{d}_k\zeta_{[2k]}.
\end{equation*}
\end{proposition}

\begin{proof}
By construction, the sequence $\mathfrak{d}_n$ is a zeon matching representation of the bipartite graph whose biadjacency matrix is $B_n$ seen in Theorem \ref{Leonardo}.  Hence, the result.
\end{proof}

\begin{example}
Let $G=(U,V,E)$ be the bipartite multigraph on $2n$ vertices $U\sqcup V$, where $|U|=|V|=n$  whose biadjacency matrix $B_n=(b_{ij})$ is determined by \begin{equation*}
b_{ij}=\begin{cases}
1&i=1, j\ne2,\\
2&i=1, j=2,\\
1&2\le i<n\,\text{\rm and }j\in\{i-1,i,i+1\},\\
1&i=n \,\text{\rm and }j\in\{n-1,n\},\\
0&\text{\rm otherwise.}
\end{cases}
\end{equation*}
Letting $\Gamma_{\mathfrak{D}}=(\delta_n)$ be the sequence defined in Proposition \ref{Leonardo zeon sequence}, the number of perfect matchings of $G$ is \begin{equation*}
\mu_G=\left< \lim_{t\to 0}\frac{d^n}{dt^n}\Wexp(t\Gamma_{\mathfrak{D}}), \zeta_{[2n]}\right>=\mathfrak{d}_n.
\end{equation*}
Graphs for $n=6, \ldots, 9$ induced by the sequence defined in \eqref{zeon Leonardo} are seen in Figure \ref{Leonardo graph sequence}.
\begin{figure}
$$\includegraphics[width=330pt]{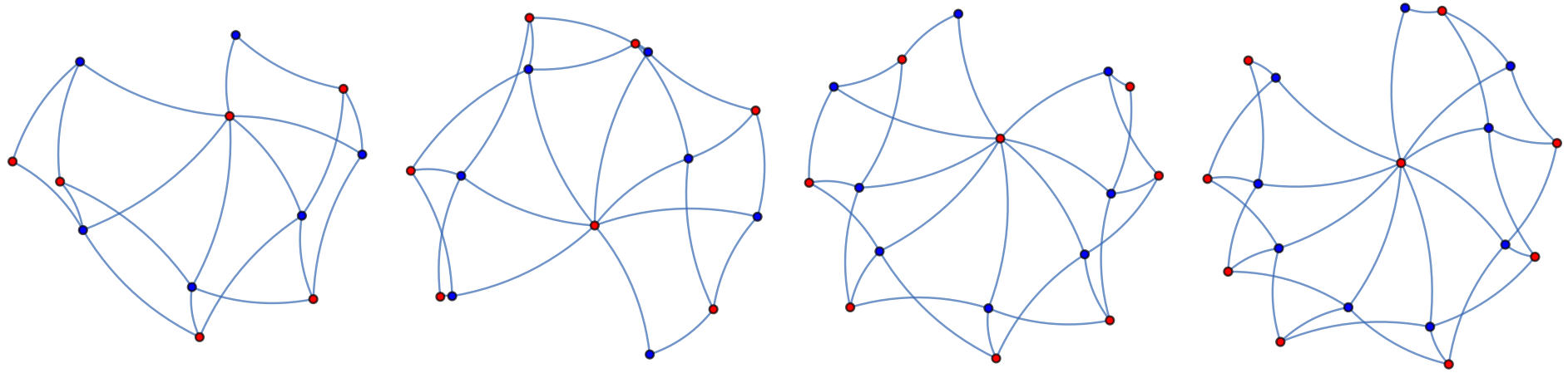}$$
\caption{Graphs associated with the zeon sequence of Proposition \ref{Leonardo zeon sequence} for $n=6, \ldots 9$.  \label{Leonardo graph sequence}}
\end{figure}
\end{example}

Finally, we arrive at a sequence of bipartite graphs whose matchings are given by a recursive sequence that is not yet found in the OEIS.  In the following proposition, we consider a sequence of bipartite multigraphs.  In each graph, a predetermined pair of vertices is incident with two edges.

\begin{proposition}\label{Aseq zeon sequence}
Let the zeon sequence $\Gamma_{\mathfrak{A}}=(\alpha_n)$ be defined recursively by \begin{align}\label{zeon aseq}
\alpha_1&=\zeta_{\{1,2\}},\nonumber\\ 
\alpha_2&=\zeta _{\{1,2\}}+2 \zeta _{\{1,4\}}+\zeta _{\{2,3\}}+\zeta _{\{3,4\}},
\nonumber\\
\alpha_3&=\zeta _{\{1,2\}}+2 \zeta _{\{1,4\}}+\zeta _{\{2,3\}}+\zeta _{\{3,4\}}+\zeta _{\{3,6\}}+\zeta _{\{4,5\}}+\zeta _{\{5,6\}},\nonumber\\
\alpha_n&=\alpha_{n-1}+\zeta_{\{2n-1,2n\}}+\zeta_{\{1,2n\}}+\zeta_{\{2n-1,2n-2\}}+\zeta_{\{2n-3,2n\}},\,\,\, (n\ge 4).
\end{align} 
Then,  \begin{equation*}
 \Wexp(t\,\Gamma_{\mathfrak{A}})=\sum_{n=0}^\infty t^n\mathfrak{a}_n\zeta_{[2n]}, 
\end{equation*}
where $(\mathfrak{a}_n)$ is the integer sequence defined by $\mathfrak{a}_0=1$, $\mathfrak{a}_1=1$, $\mathfrak{a}_2=3$, $\mathfrak{a}_3=4$, and $\mathfrak{a}_n=1+\mathfrak{a}_{n-1}+\mathfrak{a}_{n-2}$ for $n>3$.
In particular, for positive integer $k$ we have
\begin{equation*}
\lim_{t\to 0}\frac{d^k}{dt^k}\Wexp(t\,\Gamma_{\mathfrak{A}})=\mathfrak{a}_n\zeta_{[2k]}.
\end{equation*}
\end{proposition}

\begin{proof}
By construction, the sequence $\alpha_n$ is the zeon matching representation of the bipartite multigraph whose biadjacency matrix is $K_n$ seen in Theorem \ref{thm:Kn}.  Hence, the result.
\end{proof}

\begin{example}
The first five biadjacency matrices associated with the sequence $\Gamma_{\mathfrak{A}}=(\alpha_n)$ of Proposition \ref{Aseq zeon sequence} are 
\begin{align*}
\begin{pmatrix}
 1\\
\end{pmatrix},
\,\,\begin{pmatrix}
 1 & 2 \\
 1 & 1 \\
\end{pmatrix},
&\begin{pmatrix}
 1 & 2 & 0 \\
 1 & 1 & 1 \\
 0 & 1 & 1 \\
\end{pmatrix},
\begin{pmatrix}
 1 & 2 & 0 & 1 \\
 1 & 1 & 1 & 0 \\
 0 & 1 & 1 & 1 \\
 0 & 0 & 1 & 1 \\
\end{pmatrix},
\begin{pmatrix}
 1 & 2 & 0 & 1 & 1 \\
 1 & 1 & 1 & 0 & 0 \\
 0 & 1 & 1 & 1 & 0 \\
 0 & 0 & 1 & 1 & 1 \\
 0 & 0 & 0 & 1 & 1 \\
\end{pmatrix}.
\end{align*}
The graphs associated with these biadjacency matrices are see in Figure \ref{Aseq graph sequence}.  The corresponding terms of the associated zeon sequence $\Gamma_{\mathcal{A}}=(\alpha_n)$ are 
\begin{align*}
\alpha_1&= \zeta _{\{1,2\}}\\
\alpha_2&=\zeta _{\{1,2\}}+2 \zeta _{\{1,4\}}+\zeta _{\{2,3\}}+\zeta _{\{3,4\}} \\
\alpha_3&= \zeta _{\{1,2\}}+2 \zeta _{\{1,4\}}+\zeta _{\{2,3\}}+\zeta _{\{3,4\}}+\zeta _{\{3,6\}}+\zeta _{\{4,5\}}+\zeta _{\{5,6\}}\\
\alpha_4&= \zeta _{\{1,2\}}+2 \zeta _{\{1,4\}}+\zeta _{\{1,8\}}+\zeta _{\{2,3\}}+\zeta _{\{3,4\}}+\zeta _{\{3,6\}}+\zeta _{\{4,5\}}+\zeta _{\{5,6\}}\\
&\hskip10pt+\zeta _{\{5,8\}}+\zeta _{\{6,7\}}+\zeta _{\{7,8\}}\\
\alpha_5&=\zeta _{\{1,2\}}+2 \zeta _{\{1,4\}}+\zeta _{\{1,8\}}+\zeta _{\{1,10\}}+\zeta _{\{2,3\}}+\zeta _{\{3,4\}}+\zeta _{\{3,6\}}+\zeta _{\{4,5\}}\\
&\hskip10pt+\zeta _{\{5,6\}}+\zeta _{\{5,8\}}+\zeta _{\{6,7\}}
+\zeta _{\{7,8\}}+\zeta _{\{7,10\}}+\zeta _{\{8,9\}}+\zeta _{\{9,10\}}. 
\end{align*}

The first five terms of the $\Omega$-exponential $\Wexp(t\,\Gamma_{\mathfrak{A}})$ are seen here:
\begin{align*}
\Wexp(t\,\Gamma_{\mathfrak{A}})&=t \,\zeta _{\{1,2\}}+\frac{3t^2}{2}  \zeta _{\{1,2,3,4\}}+\frac{2t^3}{3}  \zeta _{\{1,2,3,4,5,6\}}+\frac{t^4}{3} \zeta _{\{1,2,3,4,5,6,7,8\}}\\
&\hskip10pt+\frac{13t^5}{120}\zeta _{\{1,2,3,4,5,6,7,8,9,10\}}+o(t^5).
\end{align*}

\begin{figure}
$$\includegraphics[width=340pt]{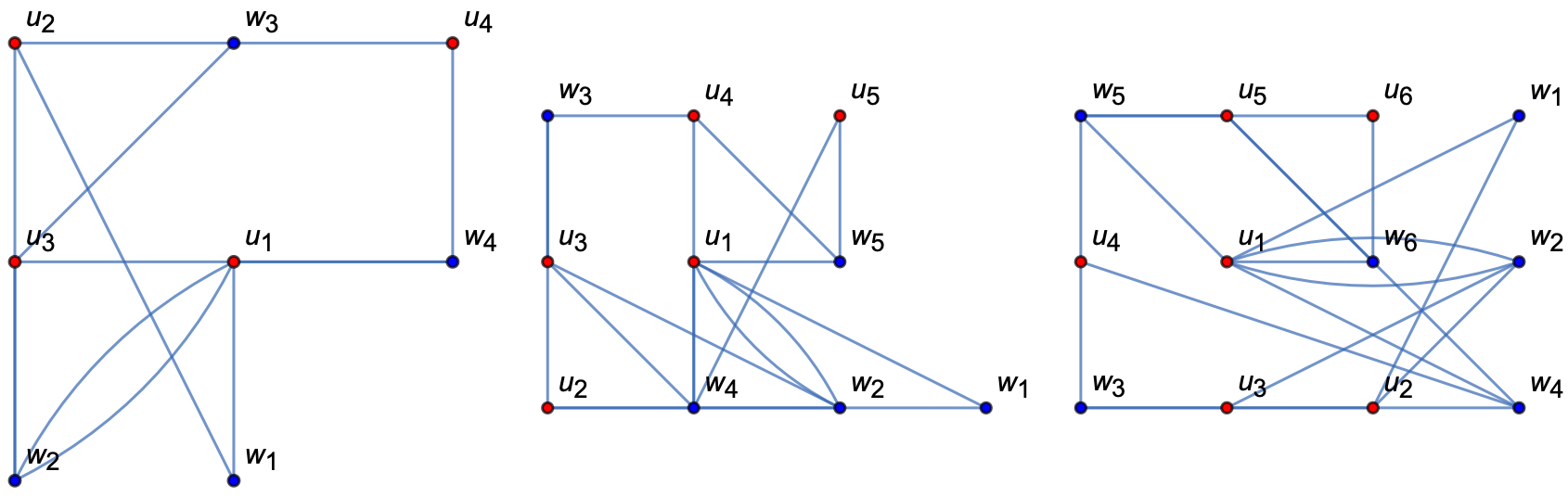}$$
\caption{Graphs associated with the zeon sequence $\Gamma_{\mathfrak{A}}=(\alpha_n)$ of Proposition \ref{Aseq zeon sequence} for $n=4, 5, 6$.  \label{Aseq graph sequence}}
\end{figure}
\end{example}

\section{Conclusion and Avenues for Future Work}\label{conclusion}

We have shown that recursively defined zeon sequences can be used to represent sequences of bipartite graphs and that their $\Omega$-exponentials reveal the numbers of perfect matchings in those graphs.  Further, the zeon sequences can be used to recover the matrix and graph sequences.

One natural extension of this work would be to better characterize those zeon sequences that generate bipartite graphs and to determine integer sequences counting the matchings of those graphs.  Generalizing to weighted bipartite graphs or bipartite multigraphs would also be worthy of investigation.  Another possible extension would be to consider zeon formulations of multipartite graphs.      

\section*{Acknowledgements}
The authors thank Arhan Canli for helpful comments on an earlier version of this paper.


\begin{thebibliography}{99}

\bibitem{armario2021}
Armario, J. ``Boolean Functions and Permanents of Sylvester Hadamard Matrices.''
\textit{Mathematics} 9, no. 2 (2021): 177.
\url{https://doi.org/10.3390/math9020177}.

\bibitem{ballot2023}
Ballot, C. J.-C., Williams, H.C..
\textit{The Lucas Sequences: Theory and Applications}.
Cham: Springer, 2023.
\url{https://doi.org/10.1007/978-3-031-37238-4}.

\bibitem{baykasoglu2014}
Baykaso\u{g}lu, A.
``A Review and Analysis of `Graph Theoretical-Matrix Permanent'
Approach to Decision Making with Example Applications.''
\textit{Artificial Intelligence Review} 42, no.~4 (2014): 573--605.
\url{https://doi.org/10.1007/s10462-012-9354-y}.

\bibitem{brualdi1991}
Brualdi, R.A., and Ryser, H.J. .
\textit{Combinatorial Matrix Theory}.
Encyclopedia of Mathematics and Its Applications, vol.~39.
Cambridge: Cambridge University Press, 1991.

\bibitem{FeinsilverMcSorleyIJC} Feinsilver, P. McSorley, J.  ``Zeons, permanents, the Johnson scheme, and generalized derangements.'' {\em Int. J. Comb.} Article ID 539030, 29 pages, (2011).\url{dx.doi.org/10.1155/2011/539030}

\bibitem{dafonseca2010}
da Fonseca, C.M.
``The $\mu$-Permanent of a Tridiagonal Matrix, Orthogonal Polynomials,
and Chain Sequences.''
\textit{Linear Algebra and its Applications} 432, no.~5 (2010):
1258--1266.
\url{https://doi.org/10.1016/j.laa.2009.10.036}.

\bibitem{koshy2011}
Koshy, T.
``Fibonacci, Lucas, and Pell Numbers, and Pascal's Triangle.''
\textit{Mathematical Spectrum} 43, no.~3 (2011): 125--132.
\url{https://oerior.uniud.it/wp-content/uploads/2018/12/Koshy2011.pdf}.

\bibitem{marcus1965}
Marcus, M., and Minc, H.
``Permanents.''
\textit{The American Mathematical Monthly} 72, no.~6
(1965): 577--591.
\url{https://doi.org/10.2307/2313846}.

\bibitem{mitchison1977}
Mitchison, G.J.
``Phyllotaxis and the Fibonacci Series.''
\textit{Science} 196, no.~4287 (1977): 270--275.
\url{https://doi.org/10.1126/science.196.4287.270}.

\bibitem{NetoJIS2014} Neto, A.F. ``Higher order derivatives of trigonometric functions, Stirling numbers of the second kind, and zeon algebra.''  {\em J. Integer Sequences}, {\bf 17} (2014), Article 14.9.3.

\bibitem{NetoJIS2015} Neto, A.F. ``Carlitz's identity for the Bernoulli numbers and zeon algebra.'' {\em  J. Integer Sequences}, {\bf 18} (2015), Article 15.5.6.


\bibitem{OEIS} OEIS Foundation Inc. (2026), The On-Line Encyclopedia of Integer Sequences, Published electronically at \url{https://oeis.org}.

\bibitem{oteles2017}
\"{O}tele\c{s}, A.
``On the Number of Perfect Matchings for Some Certain Types of
Bipartite Graphs.''
\textit{Filomat} 31, no.~15 (2017): 4809--4818.
\url{https://doi.org/10.2298/FIL1715809O}.

\bibitem{singh2017}
Singh, R., and Ravindra B. B.
``On Characteristic and Permanent Polynomials of a Matrix.''
\textit{Special Matrices} 5, no. 1 (2017): 97--112.
\url{https://doi.org/10.1515/spma-2017-0010}.

\bibitem{CAZ} Staples, G.S. {\em Clifford Algebras and Zeons: Geometry to Combinatorics and Beyond}, World Scientific Publishing, 2019. ISBN 978-981-120-257-5.

\bibitem{zeon laplacian} Staples, G.S.  ``Zeon matrix inverses and the zeon combinatorial Laplacian,''  {\em Adv. Appl. Clifford Algebras} {\bf 31}, 40 (2021).    \url{doi:10.1007/s00006-021-01152-5}.

\bibitem{vajda2008}
Vajda, S.
\textit{Fibonacci and Lucas Numbers, and the Golden Section:
Theory and Applications}.
Dover Publications, Mineola, NY. 2008.  ISBN-13: 978-0-486-46276-9.




\end{thebibliography}
\end{document}